\documentclass[11pt,a4paper]{article}
\usepackage[T1]{fontenc}
\usepackage{lmodern,amsmath,amssymb,amsthm,mathtools,geometry,booktabs,array}
\usepackage[colorlinks=true,linkcolor=blue,citecolor=blue,urlcolor=blue]{hyperref}
\newtheorem{theorem}{Theorem}[section]
\newtheorem{lemma}[theorem]{Lemma}
\newtheorem{proposition}[theorem]{Proposition}
\newtheorem{corollary}[theorem]{Corollary}
\theoremstyle{remark}
\newcommand{\R}{\mathbb R}\newcommand{\C}{\mathbb C}
\newcommand{\Om}{\Omega}\newcommand{\supp}{\operatorname{supp}}
\newcommand{\diver}{\operatorname{div}}\newcommand{\dd}{\mathrm d}
\newcommand{\N}{\mathcal N}
\newcommand{\norm}[1]{\left\|#1\right\|}
\DeclareMathOperator{\Rea}{Re}
\title{Logarithmic stability from finite-time boundary measurements\\in an infinite slab}
\author{%
  Yujian Zheng\textsuperscript{1}\qquad
  Zhiwen Duan\textsuperscript{1,2}\qquad
  Shiqi Jing\textsuperscript{1}\thanks{Corresponding author
    (\href{mailto:d202580010@hust.edu.cn}{\texttt{d202580010@hust.edu.cn}}).}\\[0.55em]
  \small\textsuperscript{1}School of Mathematics and Statistics,\\
  \small Huazhong University of Science and Technology, Wuhan 430074, Hubei, China\\[0.2em]
  \small\textsuperscript{2}Hubei Key Laboratory of Engineering Modeling and Scientific Computing,\\
  \small Huazhong University of Science and Technology, Wuhan 430074, Hubei, China}
\date{}
\begin{document}
\maketitle
\begin{abstract}
We prove logarithmic stability for an inverse convection--diffusion problem in an infinite slab from finite-time partial boundary measurements. The data consist of boundary fluxes generated by a fixed temporal profile and varying spatial Dirichlet inputs. For coefficients with fixed interior support and uniform Sobolev bounds, the data determine the exterior derivative of the drift and a gauge-invariant reaction term with a logarithmic modulus. Prescribing the divergence gives separate stability for the drift and reaction coefficients. The input and observation patches may be arbitrary nonempty open subsets of the same or opposite faces. The proof combines a quantitative passage from finite-time responses to a stationary boundary map with a Runge approximation estimate in the measured $H^{3/2}$ input norm. The first step uses sectorial semigroups on the unbounded slab; the second transfers the partial data to a full boundary map on an interior domain, where magnetic stability estimates apply. We also prove conditional stability for boundary-touching supports and uniqueness from one encoded input with unknown initial state. A tracer-transport application identifies an incompressible velocity and a reaction-loss coefficient.
\end{abstract}

\section{Introduction and main results}
Let $n\geq3$, $L>0$, and
\[
 \Om=\R^{n-1}\times(0,L),\qquad
 \Gamma_1=\{x_n=L\},\qquad\Gamma_2=\{x_n=0\}.
\]
We consider
\begin{equation}\label{eq:pde}
 \partial_tu+H_{b,c}u=0,\qquad H_{b,c}=-\Delta+b(x)\cdot\nabla+c(x),
\end{equation}
with Dirichlet input on $\Gamma_1$, homogeneous Dirichlet values on $\Gamma_2$, and observations of the normal derivative on part of either face. The coefficients are real and independent of time. Complexification is used in the analytic continuation arguments.

A motivating problem is the identification of localized transport and reaction loss by boundary tracer tests in a laterally extended layer. Advective--diffusive mass balances with reaction terms are standard in solute transport; see \cite{MOC3D}. In the idealized setting of constant known diffusivity, an incompressible transport velocity and a first-order loss coefficient lead to \eqref{eq:pde} after a change of time scale. The boundary Dirichlet input represents a prescribed concentration perturbation, and its normal derivative is proportional to the outward diffusive flux when the velocity vanishes near the faces. Corollary~\ref{cor:tracer} states the resulting identification theorem.

For autonomous parabolic equations, time analyticity connects finite-time measurements with stationary boundary data. An exact continuation argument, however, does not quantify the error in the recovered stationary map. On an unbounded domain this question also has to be treated without assuming a discrete eigenfunction expansion. We address it using a sectorial semigroup and its exponential stabilization. Our main estimate controls the stationary map by a stretched exponential of a logarithm of the inverse finite-time error. The sector geometry permits a chain of disks whose radii grow proportionally to their distance from the origin. This growth is essential to the resulting modulus.

Parabolic-to-elliptic reductions for autonomous inverse problems have a substantial history; see Canuto and Kavian \cite{CK}. Time encoding and analyticity also occur in the single-measurement results of Kian, Li, Liu and Yamamoto \cite{KLLY}. For stationary magnetic Schr\"odinger operators in a slab, Krupchyk, Lassas and Uhlmann \cite{KLU} proved partial-data recovery of the magnetic field and electric potential for observations on either the same or the opposite face. Liu and Yang \cite{LY} subsequently treated bounded magnetic potentials. The slab geometry does not reduce directly to a bounded-domain inverse problem, since Cauchy data on an artificial truncation boundary are not prescribed.

Stability of stationary partial-data problems depends on both geometry and the norm used for the data. Caro and Marinov \cite{CM} obtained stability estimates for scalar potentials in a slab using a localized solution norm for the boundary input. Partial-data magnetic stability on bounded domains was established in \cite{PR,PRT}. Our quantitative reduction is stated in a Sobolev input and normal-trace output norm. For coefficients with fixed interior support, we establish a compatible estimate of a full DN map on a bounded interior domain. The duality principle of quantitative Runge approximation, developed in \cite{RS}, is implemented with the whole-slab resolvent and the measured $H^{3/2}$ input norm. Interior Cauchy stability \cite{ARRV} controls the dual equation outside the coefficient support, where it is harmonic. The full-data theorem of Caro and Pohjola \cite{CP} applies to complex magnetic potentials and yields logarithmic coefficient stability. For supports meeting the faces we also give a nonexplicit modulus on compact classes using \cite{KLU}.

The quantitative contribution is the composition of two estimates in compatible data norms: finite-time continuation and partial-to-full boundary-map transfer. The Runge principle and the full-data magnetic estimate are taken from the existing theory; the slab realization and the norm-compatible transfer are proved here.

The first part of the paper concerns a family of experiments with one fixed time profile and varying spatial profiles. A separate result concerns a single encoded input and exact uniqueness, allowing unrelated unknown initial states. Both results use the stationary limit, but their data assumptions are different. The boundary-map continuation estimate does not require compact support. The logarithmic coefficient estimate uses fixed interior support, while the boundary-touching uniqueness result retains the geometric assumptions of the corresponding elliptic theorem.

\subsection{A quantitative boundary-map estimate}
Fix $0\leq\beta<\lambda_1:=\pi^2/L^2$.
Throughout the direct and quantitative analysis we assume
\begin{align}
 b\in W^{1,\infty}(\Om;\R^n),\qquad c\in L^\infty(\Om;\R),\label{eq:reg}\\
 c-\tfrac12\diver b\geq-\beta\quad\text{a.e. in }\Om.\label{eq:coercive}
\end{align}
Let $G\subset\Gamma_1$ and $\gamma\subset\Gamma_\sigma$, $\sigma\in\{1,2\}$, be nonempty relatively open sets. Let $X_G$ be the closure of $C_c^\infty(G)$ in $H^{3/2}(\Gamma_1)$, equipped with the inherited norm, and set $Y=L^2(\gamma)$.
For $g\in X_G$, let $P_{b,c}g\in H^2(\Om)$ solve $H_{b,c}P_{b,c}g=0$ with boundary values $g$ on $\Gamma_1$ and zero on $\Gamma_2$, and set
\[
 \N_{b,c}g=\partial_\nu P_{b,c}g|_\gamma.
\]
These maps are well defined by the direct estimates proved below.

Fix $0<a<b<t_0<t_1<T$ and a real $\eta\in C^\infty([0,\infty))$ such that $\eta=0$ on $[0,a]$ and $\eta=1$ on $[b,\infty)$. For $g\in X_G$, let $u_g$ have initial value zero and boundary input $\eta g$. Put $I=(t_0,t_1)$ and define
\[
 \mathcal R_{b,c}:X_G\longrightarrow L^2(I;Y),\qquad
 \mathcal R_{b,c}g=\partial_\nu u_g|_{I\times\gamma}.
\]
The response is recorded entirely before time $T$; the extension of the plateau beyond $T$ is used only to define the direct solution in the proof.

\begin{theorem}[Finite-time recovery of a stationary map]\label{thm:quantmain}
Fix $M>0$ and the parameters above. Let $(b_j,c_j)$, $j=1,2$, satisfy \eqref{eq:reg}--\eqref{eq:coercive} and
$\norm{b_j}_{W^{1,\infty}}+\norm{c_j}_{L^\infty}\leq M$.
There exist $C,c>0$ and $\alpha\in(0,1)$, depending only on $n,L,\beta,M,a,b,t_0,t_1,\eta$, such that
\begin{equation}\label{eq:quantmain}
 \norm{\N_{b_1,c_1}-\N_{b_2,c_2}}_{\mathcal L(X_G,Y)}
 \leq C\exp[-c\{\log(2+\varepsilon^{-1})\}^{\alpha}],
\end{equation}
where
\[
 \varepsilon=\norm{\mathcal R_{b_1,c_1}-\mathcal R_{b_2,c_2}}_{\mathcal L(X_G,L^2(I;Y))}.
\]
The right-hand side is zero at $\varepsilon=0$. No support assumption on the coefficients or relation between $G$ and $\gamma$ is required.
\end{theorem}

\subsection{Logarithmic stability for coefficients supported in the interior}
For this result the measurement sets $G\subset\Gamma_1$ and $\gamma\subset\Gamma_\sigma$ are arbitrary nonempty relatively open sets, with $\sigma\in\{1,2\}$. The unknown coefficients vanish in a fixed neighborhood of the boundary. We use the gauge invariants
\[
 (\dd b)_{jk}=\partial_jb_k-\partial_kb_j\quad(j<k),\qquad
 Q_{b,c}=c-\tfrac12\diver b+\tfrac14|b|^2.
\]
All whole-space norms below refer to zero extensions of the coefficients.

\begin{theorem}[Logarithmic coefficient stability]\label{thm:logcoeff}
Let $K\Subset\Om$ be compact, $0<s<1$, and $M>0$. Suppose that the real pairs $(b_j,c_j)$ satisfy \eqref{eq:reg}, \eqref{eq:coercive},
\begin{equation}\label{eq:logclass}
 \supp b_j\cup\supp c_j\subset K,\qquad
 \norm{b_j}_{W^{1,\infty}(\Om)}+\norm{c_j}_{L^\infty(\Om)}
 +\norm{b_j}_{H^{1+s}(\R^n)}+\norm{c_j}_{H^s(\R^n)}\leq M.
\end{equation}
Let $\varepsilon$ be the response discrepancy in Theorem~\ref{thm:quantmain}. There exist $C>0$ and $\nu>0$, depending only on the fixed geometry, $n,L,K,s,M,\beta$ and the time profile and observation interval, such that
\begin{equation}\label{eq:logcoeff}
 \norm{\dd b_1-\dd b_2}_{L^2(\Om)}+
 \norm{Q_{b_1,c_1}-Q_{b_2,c_2}}_{L^2(\Om)}
 \leq C[\log(2+\varepsilon^{-1})]^{-\nu}.
\end{equation}
The right-hand side is zero when $\varepsilon=0$. One may take
\begin{equation}\label{eq:logexponent}
 \nu=\alpha\varkappa_n\frac{s}{1+s},\qquad
 \varkappa_n=\frac{c_*}{8n^3},
\end{equation}
where $\alpha$ is furnished by Theorem~\ref{thm:quantmain} and $c_*\in(0,1)$ is a universal constant for the full-data magnetic stability estimates specified in Section~\ref{sec:explicit}.
If $\diver b_1=\diver b_2$, then also
\begin{equation}\label{eq:logfixeddiv}
 \norm{b_1-b_2}_{H^1(\Om)}+\norm{c_1-c_2}_{L^2(\Om)}
 \leq C[\log(2+\varepsilon^{-1})]^{-\nu}.
\end{equation}
\end{theorem}

\begin{corollary}[Interior H\"older classes]\label{cor:holderlog}
The conclusions of Theorem~\ref{thm:logcoeff} hold if \eqref{eq:logclass} is replaced by fixed support $K\Subset\Om$ and a uniform $C^{1,\sigma_0}(\overline\Om)\times C^{0,\sigma_0}(\overline\Om)$ bound, where $0<\sigma_0<1$. In this case any $0<s<\sigma_0$ may be used in \eqref{eq:logexponent}. For $b_1=b_2=0$, \eqref{eq:logcoeff} gives logarithmic stability of the scalar potential itself.
\end{corollary}

The proof estimates the full Dirichlet-to-Neumann map on a fixed interior domain from the measured partial map on the slab. This intermediate map is not additional data. The estimate is derived by a quantitative Runge argument using the whole-slab resolvent. It applies to the operator and its adjoint and preserves the input norm $H^{3/2}$ used in Theorem~\ref{thm:quantmain}.

\subsection{Boundary-touching supports and one encoded experiment}
For the remaining results, coefficient supports may meet the faces. Fix an open ball $B\subset\R^n$ with compact closure and relatively open sets $\gamma_{\rm in}\subset\Gamma_1$, $\gamma_{\rm out}\subset\Gamma_\sigma$ such that
\begin{equation}\label{eq:geometry}
 \overline B\cap\Gamma_1\subset\gamma_{\rm in},\qquad
 \overline B\cap\Gamma_\sigma\subset\gamma_{\rm out}.
\end{equation}
In addition to \eqref{eq:reg}--\eqref{eq:coercive}, assume
\begin{equation}\label{eq:supp}
 \supp b\cup\supp c\subset B\cap\overline\Om,\qquad b\cdot\nu=\rho\text{ on }\partial\Om,
\end{equation}
where $\rho$ is a prescribed common normal trace. Supports are relative to $\overline\Om$ and compact in the indicated ball. Coefficients may therefore meet either face. The case $\rho=0$ is included. Choose $G$ bounded with
\begin{equation}\label{eq:G}
 \overline B\cap\Gamma_1\subset G,\qquad\overline G\subset\gamma_{\rm in},
\end{equation}
and take $\gamma=\gamma_{\rm out}$. The invariants $\dd b$ and $Q_{b,c}$ are those defined above.

\begin{theorem}[Uniform conditional coefficient stability]\label{thm:compactmain}
Fix $0<\sigma_0<1$, a compact set $K\subset B\cap\overline\Om$, and $M>0$. Let $\mathcal A$ be the class satisfying \eqref{eq:coercive}, \eqref{eq:supp}, $\supp b\cup\supp c\subset K$, and
\[
 \norm b_{C^{1,\sigma_0}(\overline\Om)}+\norm c_{C^{0,\sigma_0}(\overline\Om)}\leq M.
\]
There is a bounded nondecreasing function $\omega:[0,\infty)\to[0,\infty)$, with $\omega(0)=0$ and $\omega(s)\to0$ as $s\downarrow0$, such that any pair in $\mathcal A$ satisfies
\begin{equation}\label{eq:compactmain}
 \norm{\dd b_1-\dd b_2}_{L^2(\Om)}+
 \norm{Q_{b_1,c_1}-Q_{b_2,c_2}}_{L^2(\Om)}
 \leq\omega(\Theta(\varepsilon)).
\end{equation}
Here $\varepsilon$ is as in Theorem~\ref{thm:quantmain}, and $\Theta$ denotes the right-hand side of \eqref{eq:quantmain} as a function of $\varepsilon$, with $\Theta(0)=0$. The function $\omega$ depends only on the fixed admissible class and the measurement sets; no explicit rate for $\omega$ is asserted.
If $K\Subset\Om$ and the divergence of $b$ is prescribed throughout $\mathcal A$, then also
\begin{equation}\label{eq:fixeddivstable}
 \norm{b_1-b_2}_{H^1(\Om)}+\norm{c_1-c_2}_{L^2(\Om)}
 \leq C\omega(\Theta(\varepsilon)).
\end{equation}
\end{theorem}

For a boundary function $F$ and $h\in L^2(\Om)$, write
\begin{equation}\label{eq:ibvp}
 \begin{cases}
 \partial_tu+H_{b,c}u=0 &(0,T)\times\Om,\\
 u=F &(0,T)\times\Gamma_1,\\
 u=0 &(0,T)\times\Gamma_2,\\
 u(0)=h &\Om.
 \end{cases}
\end{equation}
The normal trace for an arbitrary $L^2$ initial state is understood locally for $t>0$. We do not impose integrability at the initial time.

\begin{theorem}[One boundary experiment]\label{thm:main}
Fix $T>0$ and \eqref{eq:geometry}. There is a real input
$F\in C^\infty([0,T];H^m(\Gamma_1))$ for every integer $m\geq0$, vanishing near $t=0$ and supported spatially in a fixed compact subset of $\gamma_{\rm in}$, with the following property. For two pairs satisfying \eqref{eq:reg}, \eqref{eq:coercive}, \eqref{eq:supp}, and arbitrary $h_j\in L^2(\Om)$, let $u_j$ solve \eqref{eq:ibvp} with the same input $F$. Then
\begin{equation}\label{eq:dataeq}
 \partial_\nu u_1=\partial_\nu u_2\quad\text{on }(0,T)\times\gamma_{\rm out}
\end{equation}
implies
\begin{equation}\label{eq:conclusion}
 \dd b_1=\dd b_2,\qquad Q_{b_1,c_1}=Q_{b_2,c_2}\quad\text{in }\Om.
\end{equation}
The input depends only on the geometry and $T$, and may be multiplied by any fixed nonzero scalar.
\end{theorem}

\begin{corollary}[Prescribed divergence]\label{cor:drift}
Under Theorem~\ref{thm:main}, if the drifts have compact support in the interior of $\Om$ and $\diver b_1=\diver b_2$, then \eqref{eq:dataeq} implies $b_1=b_2$ and $c_1=c_2$.
\end{corollary}
\begin{corollary}[Scalar potential]\label{cor:scalar}
For $b=0$, the input in Theorem~\ref{thm:main} determines every real compactly supported $c\in L^\infty(\Om)$ with $c\geq-\beta$, under \eqref{eq:geometry}, even with unknown initial state.
\end{corollary}

\subsection{Application to reactive tracer transport}\label{sec:tracer}
Consider the three-dimensional layer $\Om=\R^2\times(0,L)$. Let
$\mathsf D>0$ be a known constant diffusivity, $v(x)$ a stationary transport
velocity, and $k(x)\geq0$ a first-order loss coefficient. For a dilute passive
tracer, the concentration perturbation $z(\tau,x)$ satisfies the mass balance
\begin{equation}\label{eq:tracerbalance}
 \partial_\tau z+\diver\boldsymbol J+kz=0,
 \qquad \boldsymbol J=vz-\mathsf D\nabla z.
\end{equation}
The advective, diffusive and reaction terms in this balance are standard in
solute-transport modelling; see \cite{MOC3D}. Here we adopt the idealization of
constant isotropic diffusion and assume $\diver v=0$ and
$\supp v\cup\supp k\subset K\Subset\Om$. The coefficients describe localized
transport and bulk reaction, with a quiescent, nonreactive region near the two
faces and outside a bounded lateral region. In particular, the case $v=0$
models a localized reactive inclusion in an otherwise homogeneous diffusive
layer. The infinite horizontal extent removes lateral walls from this
idealized model; no finite-width approximation estimate is assumed.

Prescribe the boundary concentration perturbation $\eta(\mathsf D\tau)g$
on the upper face and zero on the lower face, with zero initial perturbation.
The zero values represent reservoirs held at the reference concentration.
On an observation patch $\gamma$ of either face, the measured quantity is the
outward tracer flux
\begin{equation}\label{eq:tracerflux}
 \mathcal F_{v,k}g(\tau)=\boldsymbol J(\tau)\cdot\nu|_\gamma
 =-\mathsf D\,\partial_\nu z(\tau)|_\gamma,
 \qquad \tau\in I/\mathsf D.
\end{equation}
The last equality follows from the interior support assumption on $v$.
Thus the measured flux is precisely the normal observation used in this
paper, up to known units and a sign.

\begin{corollary}[Velocity and reaction loss from boundary tracer data]\label{cor:tracer}
Fix $\mathsf D>0$, $K\Subset\Om$, $0<s<1$ and $M>0$. Let the real pairs
$(v_j,k_j)$, $j=1,2$, satisfy $\diver v_j=0$, $k_j\geq0$,
$\supp v_j\cup\supp k_j\subset K$, and
\[
 \norm{v_j}_{W^{1,\infty}(\Om)}+\norm{k_j}_{L^\infty(\Om)}
 +\norm{v_j}_{H^{1+s}(\R^3)}+\norm{k_j}_{H^s(\R^3)}\leq M.
\]
For the experiments in \eqref{eq:tracerflux}, set
\[
 \varepsilon_{\rm flux}=\norm{\mathcal F_{v_1,k_1}-\mathcal F_{v_2,k_2}}
 _{\mathcal L(X_G,L^2(I/\mathsf D;Y))}.
\]
There are $C,\nu>0$, depending only on the fixed parameters and measurement
geometry, such that
\begin{equation}\label{eq:tracerstable}
 \norm{v_1-v_2}_{H^1(\Om)}+\norm{k_1-k_2}_{L^2(\Om)}
 \leq C[\log(2+\varepsilon_{\rm flux}^{-1})]^{-\nu}.
\end{equation}
The right-hand side is interpreted as zero at zero discrepancy.
\end{corollary}
\begin{proof}
Set $t=\mathsf D\tau$ and $u(t,x)=z(t/\mathsf D,x)$. Since
$\diver v_j=0$, \eqref{eq:tracerbalance} becomes \eqref{eq:pde} with
$b_j=v_j/\mathsf D$ and $c_j=k_j/\mathsf D$. These coefficients satisfy
\eqref{eq:coercive} with $\beta=0$ and have the same divergence.
Moreover,
\[
 \mathcal R_{b_j,c_j}g(t)=-\mathsf D^{-1}
 \mathcal F_{v_j,k_j}g(t/\mathsf D),\qquad
 \varepsilon=\mathsf D^{-1/2}\varepsilon_{\rm flux}.
\]
Theorem~\ref{thm:logcoeff}, followed by multiplication by $\mathsf D$,
gives \eqref{eq:tracerstable}. The known factor inside the logarithm is
absorbed into $C$.
\end{proof}

This corollary concerns identification from a family of boundary
concentration profiles with one prescribed time dependence. It gives a
conditional stability estimate for simultaneous recovery of localized
transport and reaction loss. The separate one-input theorem above concerns
exact uniqueness and does not provide the noisy-data estimate
\eqref{eq:tracerstable} from a single profile.

\subsection{Discussion of the results}
\label{sec:discussion}\label{sec:gauge}
Finite-time responses with a fixed temporal input profile control a stationary partial boundary map with an explicit modulus. For coefficients supported away from the faces, quantitative approximation of local solutions then yields logarithmic stability for the gauge invariants, and for the drift and reaction separately when the divergence is prescribed. In the tracer model, incompressibility provides this divergence condition and the measured normal trace is the physical outward flux after a known rescaling.

The stability estimate compares boundary response operators over a family of spatial profiles. A finite collection of actuators or sensors, errors caused by lateral truncation, and a reconstruction algorithm require additional analysis. Supports meeting the boundary are covered here by conditional stability without a specified rate and by the separate single-input uniqueness result. Extending the explicit coefficient estimate to such supports remains a further question.

Let $\psi\in C_c^\infty(\Om;\R)$. Direct differentiation gives
\begin{equation}\label{eq:gauge}
 e^{-\psi}H_{b,c}e^\psi=H_{\widetilde b,\widetilde c},\qquad
 \widetilde b=b-2\nabla\psi,\quad
 \widetilde c=c+b\cdot\nabla\psi-\Delta\psi-|\nabla\psi|^2.
\end{equation}
Since $\psi$ is independent of time,
\[
 e^{-\psi}(\partial_t+H_{b,c})e^\psi=\partial_t+H_{\widetilde b,\widetilde c}.
\]
Near the boundary $\psi=0$, so Dirichlet and ordinary Neumann values are unchanged. The initial state becomes $e^{-\psi}h$; thus full boundary data cannot distinguish these coefficient pairs on an unrestricted class.

Substitution into the definition of $Q_{b,c}$ gives
\[
 \dd\widetilde b=\dd b,\qquad Q_{\widetilde b,\widetilde c}=Q_{b,c}.
\]
This ambiguity also occurs in coercive classes. Starting from $b=c=0$, choose a nonzero compactly supported $\psi$ with sufficiently small gradient. Then
\[
 \widetilde c-\tfrac12\diver\widetilde b=-|\nabla\psi|^2\geq-\beta.
\]
Fixing the divergence removes this ambiguity under the support assumption of Corollary~\ref{cor:drift}.

The proof is organized as follows. Section~\ref{sec:direct} develops the direct problem, stationary states and time analyticity. Section~\ref{sec:quant} establishes the quantitative passage from finite-time measurements to the stationary map. Section~\ref{sec:explicit} proves the logarithmic coefficient estimates. Section~\ref{sec:uniqueness} treats the encoded experiment, the elliptic uniqueness argument, compact-class stability and the magnetic heat equation.

\section{The direct problem and time analyticity}\label{sec:direct}
All Hilbert inner products are linear in their first argument. We write $X=L^2(\Om;\C)$ and $V=H_0^1(\Om;\C)$. The Poincar\'e inequality on each vertical interval gives
\begin{equation}\label{eq:poincare}
 \lambda_1\norm{v}_{X}^{2}\leq\norm{\partial_nv}_{X}^{2}
 \leq\norm{\nabla v}_{X}^{2},\qquad v\in V.
\end{equation}
It follows first for smooth compactly supported functions by the one-dimensional Dirichlet inequality and Fubini, and then for $V$ by density.

\subsection{Sectorial realization and regularity}
Define the form
\begin{equation}\label{eq:form}
 a(v,w)=\int_\Om\big(\nabla v\cdot\nabla\overline w+
 b\cdot\nabla v\,\overline w+c\,v\overline w\big)\,\dd x,
 \qquad v,w\in V.
\end{equation}
There is no magnetic conjugation in this definition. Integration by parts in the drift term yields
\begin{equation}\label{eq:energy}
 \Rea a(v,v)=\norm{\nabla v}_X^2+
 \int_\Om(c-\tfrac12\diver b)|v|^2\,\dd x
 \geq(1-\beta/\lambda_1)\norm{\nabla v}_X^2.
\end{equation}
The integration by parts follows by density from $C_c^\infty(\Om)$, using the $W^{1,\infty}$ bound on $b$. In particular,
\[
 \Rea a(v,v)\geq(\lambda_1-\beta)\norm v_X^2.
\]

\begin{lemma}\label{lem:operator}
The form \eqref{eq:form} defines a closed densely defined operator $H$ in $X$ with
\[
 D(H)=H^2(\Om)\cap H_0^1(\Om),\qquad Hv=H_{b,c}v.
\]
It is invertible, and $-H$ generates an analytic semigroup $S(t)$ such that
\begin{equation}\label{eq:sgbounds}
 \norm{S(t)}_{X\to X}\leq e^{-(\lambda_1-\beta) t},\qquad
 \norm{HS(t)}_{X\to X}\leq C(1+t^{-1})e^{-(\lambda_1-\beta) t/2}.
\end{equation}
Moreover,
\begin{equation}\label{eq:graph}
 \norm v_{H^2(\Om)}\leq C\big(\norm{Hv}_X+\norm v_X\big),\qquad v\in D(H).
\end{equation}
The adjoint $H^*$ has the same domain and is the Dirichlet realization of
$-\Delta-b\cdot\nabla+c-\diver b$.
\end{lemma}
\begin{proof}
Continuity and \eqref{eq:energy} give a closed sectorial form on $V$, since
\[
 |\Im a(v,v)|\leq\norm b_\infty\norm{\nabla v}_X\norm v_X
 \leq C\Rea a(v,v),\qquad
 \Rea a(v,v)\asymp\norm v_{H^1}^2.
\]
The representation and analytic-semigroup theorems \cite{Ouhabaz}, together with Lax--Milgram, yield
\[
 H^{-1}:X\longrightarrow V,\qquad S(t)=e^{-tH},\qquad
 \norm{H^{-1}}_{X\to X}\leq(\lambda_1-\beta)^{-1}.
\]
The energy inequality implies
\[
 \frac{\dd}{\dd t}\norm{S(t)v}_X^2
 =-2\Rea a(S(t)v,S(t)v)\leq-2(\lambda_1-\beta)\norm{S(t)v}_X^2,
\]
and hence the first bound in \eqref{eq:sgbounds}, initially on $D(H)$ and then by density.

If $v\in D(H)$ and $Hv=f$, then
\[
 -\Delta v=f-b\cdot\nabla v-cv\in X.
\]
Fourier transformation in $x'$ and the Dirichlet sine expansion in $x_n$ give
\[
 \norm v_{H^2}\leq C(\norm{\Delta v}_X+\norm v_X),\qquad
 |\xi'|^2(k\pi/L)^2\leq\tfrac14\big(|\xi'|^2+(k\pi/L)^2\big)^2.
\]
The second inequality controls the mixed derivatives. Integration by parts proves the converse inclusion, so
\[
 D(H)=H^2(\Om)\cap V.
\]
Absorb the first-order term using
\[
 \norm{\nabla v}_X\leq\epsilon\norm v_{H^2}+C_\epsilon\norm v_X
\]
to obtain \eqref{eq:graph}. Analyticity and the semigroup property give
\[
 \norm{HS(t)}_{X\to X}\leq C/t\quad(0<t\leq1),\qquad
 HS(t)=HS(1)S(t-1)\quad(t>1).
\]
These imply the second bound in \eqref{eq:sgbounds}. Finally, integration by parts yields
\[
 H^*=-\Delta-b\cdot\nabla+c-\diver b,\qquad
 (c-\diver b)-\tfrac12\diver(-b)=c-\tfrac12\diver b.
\]
The same coercivity and regularity arguments therefore give $D(H^*)=H^2(\Om)\cap V$.
\end{proof}

The trace theorem and \eqref{eq:graph} give a bounded map
\begin{equation}\label{eq:C}
 C_\gamma:D(H)\longrightarrow L^2(\gamma),\qquad
 C_\gamma v=\partial_\nu v|_\gamma,
\end{equation}
where $D(H)$ has its graph norm and $\gamma$ is any relatively open portion of either face. The normal trace is in fact in $H^{1/2}$ on the full flat boundary. For a nonzero Dirichlet solution we use the same normal trace on $H^2(\Om)$.

\subsection{Lifting the boundary values}
For $g\in H^{3/2}(\Gamma_1)$ there is a bounded extension $Eg\in H^2(\Om)$ with traces $g$ on $\Gamma_1$ and zero on $\Gamma_2$. For completeness, write $y=L-x_n$, choose a smooth $\chi(y)$ equal to one near zero and vanishing near $L$, and set
\[
 \widehat{Eg}(\xi',y)=\chi(y)e^{-y\langle\xi'\rangle}\widehat g(\xi'),
 \qquad \langle\xi'\rangle=(1+|\xi'|^2)^{1/2}.
\]
Integrating the squares of all derivatives up to order two in $y$ gives
$\norm{Eg}_{H^2}\leq C\norm g_{H^{3/2}}$; the highest power is
$\langle\xi'\rangle^3$. The cutoff terms have no larger order. This also supplies a bounded lifting for time-dependent data.

For $f\in C^2([0,T];H^{3/2}(\Gamma_1))$ vanishing near zero, let $z=Ef$ and put
\begin{equation}\label{eq:mild}
 u(t)=z(t)+S(t)h+\int_0^t S(t-s)G(s)\,\dd s,
 \qquad G=-z'-H_{b,c}z.
\end{equation}
Here $H_{b,c}z$ is the differential expression, not the Dirichlet operator applied to a vector outside its domain. The integral is a Bochner integral in $X$.

\begin{proposition}\label{prop:wellposed}
Formula \eqref{eq:mild} gives the unique energy solution of \eqref{eq:ibvp}. For every $\varepsilon>0$ it belongs to
$C([\varepsilon,T];H^2(\Om))\cap C^1([\varepsilon,T];X)$, and its boundary flux is well defined there. The solution depends continuously on $h\in X$ and on $f$ in the above $C^2$ norm.
\end{proposition}
\begin{proof}
Subtracting the lifting gives
\[
 v=u-z,\qquad v'+Hv=G,\qquad v(0)=h,
\]
so variation of constants yields \eqref{eq:mild} and continuity in $X$. To control the graph norm, write
\begin{align*}
 H\int_0^tS(t-s)G(s)\,\dd s
 &=\int_0^tHS(t-s)(G(s)-G(t))\,\dd s\\
 &\quad +(I-S(t))G(t).
\end{align*}
Approximation from smooth $D(H)$-valued data proves the identity. Since $G\in C^1([0,T];X)$,
\[
 \norm{HS(t-s)(G(s)-G(t))}_X\leq C_T\norm{G'}_{C([0,T];X)}.
\]
Dominated convergence and \eqref{eq:graph} now imply, for every $\varepsilon>0$,
\[
 u\in C([\varepsilon,T];H^2),\qquad
 u'=-H_{b,c}u\in C([\varepsilon,T];X).
\]
The same bounds give continuous dependence and local energy regularity. For a difference $v$ with zero data,
\[
 \tfrac12\frac{\dd}{\dd t}\norm v_X^2+(\lambda_1-\beta)\norm v_X^2\leq0,
 \qquad v(0)=0,
\]
which proves uniqueness.
\end{proof}

\subsection{Stationary states and time analyticity}
For $g\in H^{3/2}(\Gamma_1)$ define
\begin{equation}\label{eq:poisson}
 P g=Eg-H^{-1}(H_{b,c}Eg).
\end{equation}
This formula is independent of the lifting, belongs to $H^2(\Om)$, and satisfies
\begin{equation}\label{eq:stationary}
 H_{b,c}Pg=0,\qquad Pg|_{\Gamma_1}=g,\qquad Pg|_{\Gamma_2}=0.
\end{equation}
Indeed the second term in \eqref{eq:poisson} has zero Dirichlet trace, and applying the differential expression cancels the two terms. Uniqueness follows from invertibility of $H$. The elliptic partial DN map at energy zero is
\[
 \N_{b,c}^{\gamma}g=\partial_\nu Pg|_\gamma.
\]
It is bounded from $H^{3/2}(\Gamma_1)$ to $L^2(\gamma)$.

\begin{lemma}[A ramp followed by a plateau]\label{lem:plateau}
Let $0<a<b$ and let $\eta\in C^\infty([0,\infty))$ vanish for $t\leq a$ and equal one for $t\geq b$. Let $u_g$ solve \eqref{eq:ibvp} with zero initial value and boundary input $\eta(t)g$. Then, for $t>b$,
\begin{equation}\label{eq:plateau}
 u_g(t)=Pg+S(t-b)(u_g(b)-Pg).
\end{equation}
Consequently $t\mapsto\partial_\nu u_g(t)|_\gamma$ is real analytic with values in $L^2(\gamma)$ on $(b,\infty)$, and
\begin{equation}\label{eq:limit}
 \partial_\nu u_g(t)|_\gamma\longrightarrow\N_{b,c}^{\gamma}g
 \quad\text{in }L^2(\gamma)\quad(t\to\infty).
\end{equation}
The same statements hold with an arbitrary initial state.
\end{lemma}
\begin{proof}
For $t\geq b$, the difference $u_g-Pg$ has zero Dirichlet trace and solves
\[
 (\partial_t+H)(u_g-Pg)=0,\qquad
 (u_g-Pg)(b)=u_g(b)-Pg\in X.
\]
Semigroup uniqueness proves \eqref{eq:plateau}. For a fixed $t_0>b$, factor
\[
 S(t-b)=S((t_0-b)/2)S(t-b-(t_0-b)/2)
\]
near $t_0$. The first factor maps $X$ boundedly into $D(H)$, while the second is analytic in $\mathcal L(X)$.
Consequently,
\[
 t\longmapsto C_\gamma S(t-b)(u_g(b)-Pg)
\]
is analytic into $L^2(\gamma)$. Moreover, \eqref{eq:sgbounds} gives
\[
 \norm{C_\gamma S(t-b)(u_g(b)-Pg)}_Y
 \leq C(1+(t-b)^{-1})e^{-(\lambda_1-\beta)(t-b)/2}\norm{u_g(b)-Pg}_X
 \longrightarrow0.
\]
This proves \eqref{eq:limit}; an arbitrary initial state changes only $u_g(b)-Pg$.
\end{proof}

\begin{lemma}[Equality on a time interval]\label{lem:analytic}
For two admissible operators, let $u_{j,g}$ denote their zero-initial-state responses to the same ramp $\eta g$. If their fluxes agree on $J\times\gamma$ for a nonempty open interval $J\subset(b,\infty)$, then
\[
 \N_{b_1,c_1}^{\gamma}g=\N_{b_2,c_2}^{\gamma}g.
\]
Moreover the two flux responses agree for every $t>b$.
\end{lemma}
\begin{proof}
For every $\ell\in Y^*$, the scalar analytic function
\[
 t\longmapsto\ell\big(\partial_\nu u_{1,g}(t)|_\gamma
                   -\partial_\nu u_{2,g}(t)|_\gamma\big)
\]
vanishes on $J$, hence on $(b,\infty)$ by the identity theorem. Duality gives equality of the fluxes there, and \eqref{eq:limit} gives equality of the stationary maps on $g$.
\end{proof}

\section{Quantitative continuation to the stationary map}\label{sec:quant}
We first establish a continuation estimate in a Banach space. For $0<\vartheta<\pi/2$ and $r_*>0$, put
\[
 \mathcal S_{\vartheta,r_*}=\{z\in\C:|\arg z|<\vartheta,\ |z|>r_*\}.
\]
The constants in the following lemma are independent of the Banach space.

\begin{lemma}[Continuation of a stabilizing holomorphic function]\label{lem:quantanalytic}
Let $Y$ be a complex Banach space, let $J=(s_0,s_1)$ with $r_*<s_0<s_1$, and let $\lambda>0$. Suppose $F:\mathcal S_{\vartheta,r_*}\to Y$ is holomorphic and, for some $M_0>0$ and $F_\infty\in Y$,
\begin{equation}\label{eq:abstractbounds}
 \sup_{z\in\mathcal S_{\vartheta,r_*}}\norm{F(z)}_Y\leq M_0,
 \qquad \norm{F(x)-F_\infty}_Y\leq M_0e^{-\lambda x}\quad(x>s_0).
\end{equation}
There are $C,c>0$ and $\alpha\in(0,1)$, depending only on $\vartheta,r_*,J,\lambda$, such that, with $\varepsilon=\norm F_{L^2(J;Y)}$,
\begin{equation}\label{eq:abstractmodulus}
 \norm{F_\infty}_Y\leq CM_0
 \exp\left[-c\left\{\log\left(2+\frac{M_0}{\varepsilon}\right)\right\}^{\alpha}\right].
\end{equation}
The right-hand side is interpreted as zero when $\varepsilon=0$.
\end{lemma}
\begin{proof}
We may divide $F$ and $F_\infty$ by $M_0$ and assume $M_0=1$.
Choose $x_0\in J$ and $0<\eta<1/8$ sufficiently small that, with $r_0=\eta x_0$,
\[
 [x_0-2r_0,x_0+2r_0]\subset J,\qquad
 \overline{D(x,4\eta x)}\subset\mathcal S_{\vartheta,r_*}
 \quad(x\geq x_0).
\]
It suffices to choose $4\eta<\sin\vartheta$ and $(1-4\eta)x_0>r_*$, independently of $F$. Cauchy's estimate gives
\[
 \sup_{|x-x_0|\leq2r_0}\norm{F'(x)}_Y\leq C/r_0.
\]
For $|x-x_0|\leq r_0$, the reverse triangle inequality therefore yields
\[
 \norm{F(y)}_Y\geq\tfrac12\norm{F(x)}_Y
 \quad\text{if }|y-x|\leq cr_0\norm{F(x)}_Y.
\]
Choose $c$ so that this interval stays in the larger segment. Integration gives
\[
 \varepsilon^2\geq cr_0\norm{F(x)}_Y^3,
\]
and consequently
\begin{equation}\label{eq:smallsegment}
 m:=\sup_{|x-x_0|\leq r_0}\norm{F(x)}_Y\leq C\varepsilon^{2/3}.
\end{equation}

Let $E=[x_0-r_0,x_0+r_0]$. The harmonic measure $\omega$ of this slit satisfies
\[
 \Delta\omega=0\text{ in }D(x_0,4r_0)\setminus E,\qquad
 \omega|_E=1,\qquad\omega|_{\partial D(x_0,4r_0)}=0.
\]
Rescaling, the minimum principle and regularity at the slit endpoints give a fixed $\theta_0\in(0,1)$ with
\[
 \omega\geq\theta_0\quad\text{on }\overline{D(x_0,r_0)}\setminus E.
\]
For $\ell\in Y^*$, $\norm\ell\leq1$, apply the subharmonic maximum principle to $\log|\ell(F)|$, regularized at its zeros. The two-constants inequality gives
\[
 |\ell(F(z))|\leq m^{\omega(z)}\leq m^{\theta_0}
 \quad (z\in D(x_0,r_0)\setminus E),
\]
provided $m<1$. On $E$ the same bound follows directly. The dual characterization of the norm therefore yields
\begin{equation}\label{eq:smalldisk}
 \sup_{D(x_0,r_0)}\norm F_Y\leq m^{\theta_0}.
\end{equation}
If $m=0$, scalarization and the identity theorem give $F=0$ everywhere and the conclusion follows. We henceforth suppose $\varepsilon>0$ is small enough that $0<m<1$.

Set $a=1+\eta/2$, $x_k=x_0a^k$ and $r_k=\eta x_k$. Then
\[
 x_{k+1}-x_k=r_k/2,\qquad
 D(x_{k+1},r_{k+1})\subset D(x_k,2r_k).
\]
The outer disks $D(x_k,4r_k)$ all lie in the sector. Hadamard's three-circles theorem, applied to $\ell(F)$ and then taking the supremum over $\ell$, implies
\[
 \sup_{D(x_k,2r_k)}\norm F_Y
 \leq\left(\sup_{D(x_k,r_k)}\norm F_Y\right)^{1/2}.
\]
Induction with \eqref{eq:smalldisk} gives
\begin{equation}\label{eq:chain}
 \norm{F(x_k)}_Y\leq m^{\theta_0 2^{-k}}
 \leq \exp[-c2^{-k}\log(1/\varepsilon)]
\end{equation}
for all $k\geq0$, after reducing the fixed upper threshold for $\varepsilon$ in \eqref{eq:smallsegment}.

Set $\alpha=\log a/\log(2a)\in(0,1)$. For sufficiently small $\varepsilon>0$, choose $k\geq0$ with
\[
 a^k\leq[\log(1/\varepsilon)]^\alpha<a^{k+1}.
\]
Then
\[
 x_k\geq(x_0/a)[\log(1/\varepsilon)]^\alpha,\qquad
 2^{-k}\geq[\log(1/\varepsilon)]^{\alpha-1}.
\]
Equations \eqref{eq:abstractbounds} and \eqref{eq:chain} yield
\[
 \norm{F_\infty}_Y\leq\norm{F(x_k)}_Y+e^{-\lambda x_k}
 \leq2\exp\big(-c[\log(1/\varepsilon)]^\alpha\big).
\]
This is \eqref{eq:abstractmodulus} for small $\varepsilon$. For the remaining values, enlarge $C$ and use $\norm{F_\infty}_Y\leq1$. When $\varepsilon=0$, continuity on $J$ followed by the identity theorem proves the zero-error assertion.
\end{proof}

\begin{lemma}[Uniform sector estimates]\label{lem:uniform}
Assume \eqref{eq:reg}, \eqref{eq:coercive} and
\begin{equation}\label{eq:uniformbound}
 \norm b_{W^{1,\infty}(\Om)}+\norm c_{L^\infty(\Om)}\leq M.
\end{equation}
There exist $\vartheta\in(0,\pi/2)$ and $C>0$, depending only on $n,L,\beta,M$, such that the semigroups extend holomorphically to $|\arg z|<\vartheta$ and satisfy
\begin{equation}\label{eq:complexsg}
 \norm{S(z)}_{L^2\to L^2}+\norm{HS(z)}_{L^2\to L^2}
 \leq C(1+|z|^{-1})e^{-(\lambda_1-\beta)\Rea z/2}.
\end{equation}
The graph estimate \eqref{eq:graph} and $\norm P_{H^{3/2}\to H^2}\leq C$ hold uniformly on this class.
\end{lemma}
\begin{proof}
By \eqref{eq:energy}, the shifted form satisfies
\begin{align*}
 \Rea\big(a(v,v)-\tfrac{\lambda_1-\beta}{2}\norm v_X^2\big)
 &\geq\tfrac12(1-\beta/\lambda_1)\norm{\nabla v}_X^2,\\
 |\Im a(v,v)|&\leq M\lambda_1^{-1/2}\norm{\nabla v}_X^2.
\end{align*}
Thus $H-(\lambda_1-\beta)/2$ has a uniform sectorial angle. On a smaller common sector, the form theorem and Cauchy's estimate give
\[
 \norm{e^{-z(H-(\lambda_1-\beta)/2)}}_{X\to X}\leq C,\qquad
 \norm{(H-(\lambda_1-\beta)/2)e^{-z(H-(\lambda_1-\beta)/2)}}_{X\to X}\leq C|z|^{-1}.
\]
Multiplication by $e^{-(\lambda_1-\beta) z/2}$ proves \eqref{eq:complexsg}. The absorption argument above is uniform, hence
\[
 \norm v_{H^2}\leq C(\norm{Hv}_X+\norm v_X),\qquad
 \norm{H^{-1}}_{X\to H^2}\leq C.
\]
The fixed lifting and \eqref{eq:poisson} then yield
\[
 \norm{Pg}_{H^2}\leq\norm{Eg}_{H^2}
       +C\norm{H_{b,c}Eg}_X\leq C\norm g_{H^{3/2}}.
\]
\end{proof}

\begin{proof}[Proof of Theorem~\ref{thm:quantmain}]
Write $C_jv=\partial_\nu v|_{\gamma}$ and $P_j$ for the two stationary Poisson operators. For the response to $\eta g$, the equation with stationary lifting gives
\[
 u_{j,g}(t)=\eta(t)P_jg-\int_0^tS_j(t-s)\eta'(s)P_jg\,\dd s.
\]
In particular, with $W_jg=u_{j,g}(b)-P_jg$,
\begin{equation}\label{eq:W}
 W_jg=-\int_0^bS_j(b-s)\eta'(s)P_jg\,\dd s,
 \qquad \norm{W_jg}_{L^2}\leq C\norm g_{H^{3/2}}.
\end{equation}
For $\norm g_{X_G}=1$, define the $Y$-valued function
\begin{equation}\label{eq:Fg}
 F_g(z)=(\N_1-\N_2)g+C_1S_1(z)W_1g-C_2S_2(z)W_2g.
\end{equation}
It is holomorphic on the common sector away from its vertex. Choose $0<r_*<t_0-b$. Lemma~\ref{lem:uniform}, the trace bound, and \eqref{eq:W} imply, uniformly in $g$ and the coefficients,
\[
 \sup_{\mathcal S_{\vartheta,r_*}}\norm{F_g(z)}_Y\leq C,
 \qquad \norm{F_g(x)-(\N_1-\N_2)g}_Y\leq Ce^{-(\lambda_1-\beta) x/2}
 \quad(x>t_0-b).
\]
For real $z=t-b$, Lemma~\ref{lem:plateau} identifies $F_g(t-b)$ with $(\mathcal R_1-\mathcal R_2)g(t)$. Hence
\[
 \norm{F_g}_{L^2((t_0-b,t_1-b);Y)}\leq\varepsilon.
\]
Lemma~\ref{lem:quantanalytic} now proves \eqref{eq:quantmain} for each unit vector $g$, with the same constants. Taking the supremum over $g$ proves the operator estimate. This last step uses the operator norm into $L^2(I;Y)$; it does not interchange a spatial supremum with a time integral.
\end{proof}

\begin{corollary}[Subtraction of a repeatable initial state]\label{cor:baseline}
For each fixed operator, suppose the experiments with boundary inputs $\eta g$ and zero have the same initial state $h\in L^2(\Om)$. The difference of their measured fluxes is $\mathcal R_{b,c}g$. Theorem~\ref{thm:quantmain} therefore applies to these differences without knowledge of $h$.
\end{corollary}
\begin{proof}
Subtract the two direct equations. The difference has zero initial state and boundary input $\eta g$, so uniqueness of the direct problem identifies it with $u_g$.
\end{proof}

\begin{proposition}[Transfer of elliptic stability]\label{prop:transfer}
Let $\mathcal A$ be a coefficient class satisfying the hypotheses of Theorem~\ref{thm:quantmain}. Let $d$ be a nonnegative measure of coefficient discrepancy and suppose that, for all pairs in $\mathcal A$,
\begin{equation}\label{eq:ellipticinterface}
 d((b_1,c_1),(b_2,c_2))\leq
 \omega_E\big(\norm{\N_1-\N_2}_{\mathcal L(X_G,Y)}\big),
\end{equation}
where $\omega_E$ is nondecreasing and vanishes at zero. Then
\[
 d((b_1,c_1),(b_2,c_2))\leq\omega_E(\Theta(\varepsilon)),
 \qquad \Theta(s)=C\exp[-c\{\log(2+s^{-1})\}^{\alpha}],\quad\Theta(0)=0.
\]
In particular, a bound $\omega_E(s)\leq C[\log(2+s^{-1})]^{-\mu}$ for small $s>0$ yields
\[
 d((b_1,c_1),(b_2,c_2))\leq C[\log(2+\varepsilon^{-1})]^{-\alpha\mu}.
\]
A bound $\omega_E(s)\leq C[\log\log(e^e+s^{-1})]^{-\mu}$ yields
\[
 d((b_1,c_1),(b_2,c_2))\leq
 C[\log\log(e^e+\varepsilon^{-1})]^{-\mu}.
\]
These implications require \eqref{eq:ellipticinterface} in the indicated spaces $X_G$ and $Y$.
\end{proposition}
\begin{proof}
Substitution of \eqref{eq:quantmain} into \eqref{eq:ellipticinterface} gives the first assertion. For sufficiently small $s>0$,
\[
 \log(2+\Theta(s)^{-1})\geq c[\log(2+s^{-1})]^\alpha,
\]
and therefore
\[
 \log\log(e^e+\Theta(s)^{-1})\geq c\log\log(e^e+s^{-1}).
\]
Applying the assumed elliptic modulus proves the two rates. Enlarging the constants covers any fixed bounded interval of errors.
\end{proof}

\section{Explicit coefficient stability in the interior}\label{sec:explicit}
In this section all coefficient supports lie in a fixed compact set $K\Subset\Om$. Choose smooth bounded domains
\begin{equation}\label{eq:nested}
 K\Subset D_0\Subset D_1\Subset D_2\Subset D\Subset\Om
\end{equation}
such that $\Om\setminus\overline{D_1}$ is connected. These domains can be chosen convex: the convex hull of $K$ remains at positive distance from both faces. All auxiliary sets below are fixed independently of the coefficients. Constants are uniform under the $W^{1,\infty}\times L^\infty$ bound and the coercivity assumption.

\subsection{A boundary-to-interior estimate for the dual equation}
Let $\Gamma\subset\Gamma_\ell$ be a bounded open disk with closure in one of the available boundary patches. Let $\mathcal X_\Gamma$ be the closure of $C_c^\infty(\Gamma)$ in $H^{3/2}(\Gamma_\ell)$, equipped with the inherited norm.
The Poisson operator $P_H$ now allows data on either face and zero values on the other; the direct estimates follow from Section~2, with reflection in $x_n=L/2$ when necessary. Define
\[
 \mathcal A_H:\mathcal X_\Gamma\longrightarrow L^2(D_1),\qquad
 \mathcal A_H f=(P_H f)|_{D_1}.
\]
For $r\in L^2(D_1)$, extend $r$ by zero and set
\[
 w=(H^*)^{-1}(\mathbf1_{D_1}r).
\]
The star on $\mathcal A_H$ denotes the Hilbert space adjoint with the inherited $H^{3/2}$ inner product. The equality of its norm with the norm of the corresponding dual functional will be used below.

\begin{lemma}\label{lem:dualholder}
Let $S\Subset\Om\setminus\overline{D_1}$ be compact. There are $C>0$ and $\theta\in(0,1)$ such that
\begin{equation}\label{eq:dualholder}
 \norm w_{L^2(S)}\leq C\norm r_{L^2(D_1)}^{1-\theta}\norm{\mathcal A_H^*r}_{\mathcal X_\Gamma}^{\theta}.
\end{equation}
The constants depend on $S$, the fixed sets, and the uniform direct-problem bounds. The estimate also holds with $H$ replaced by $H^*$.
\end{lemma}
\begin{proof}
The global graph estimate gives
\begin{equation}\label{eq:dualglobal}
 \norm w_{H^2(\Om)}\leq C\norm r_{L^2(D_1)}.
\end{equation}
Green's identity, followed by density, gives
\[
 (\mathcal A_H f,r)_{L^2(D_1)}
 =-\int_\Gamma f\,\overline{\partial_\nu w}\,\dd S.
\]
The drift boundary term vanishes because $w|_{\partial\Om}=0$. Hence
\begin{equation}\label{eq:dualtrace}
 \left|\int_\Gamma f\,\overline{\partial_\nu w}\,\dd S\right|
 \leq\norm{\mathcal A_H^*r}_{\mathcal X_\Gamma}\norm f_{H^{3/2}(\Gamma_\ell)}.
\end{equation}
Choose $\Gamma_0\Subset\Gamma$ and $\zeta\in C_c^\infty(\Gamma)$ equal to one near $\overline{\Gamma_0}$. Testing with $\zeta\varphi$ and applying the trace theorem yield
\[
 \norm{\zeta\partial_\nu w}_{H^{-3/2}(\Gamma_\ell)}
 \leq C\norm{\mathcal A_H^*r}_{\mathcal X_\Gamma},\qquad
 \norm{\zeta\partial_\nu w}_{H^{1/2}(\Gamma_\ell)}\leq C\norm r_{L^2(D_1)}.
\]
Interpolation on the flat hyperplane gives
\begin{equation}\label{eq:traceinterpolation}
 \norm{\partial_\nu w}_{H^{-1/2}(\Gamma_0)}
 \leq C\norm r_{L^2(D_1)}^{1/2}\norm{\mathcal A_H^*r}_{\mathcal X_\Gamma}^{1/2}.
\end{equation}
Here $H^{-1/2}(\Gamma_0)$ is dual to the supported $H^{1/2}$ space with its inherited hyperplane norm. Outside $\overline{D_1}$,
\[
 -\Delta w=0,
\]
because the source and coefficients vanish there. Connectedness permits a fixed bounded Lipschitz domain $V$ with
\[
 S\Subset V\subset\Om\setminus\overline{D_1},\qquad
 \Gamma_0\subset\partial V,\qquad w|_{\Gamma_0}=0.
\]
One constructs $V$ by joining an inward collar of $\Gamma_0$ to a neighborhood of $S$ by finitely many interior tubes.

The harmonic case of \cite[Theorem~1.7]{ARRV}, applied to a compact connecting subregion of $V$, gives
\begin{equation}\label{eq:externalcauchy}
 \norm w_{L^2(S)}\leq
 C\big(\norm w_{H^1(V)}+\norm{\partial_\nu w}_{H^{-1/2}(\Gamma_0)}\big)^{1-\vartheta}
 \norm{\partial_\nu w}_{H^{-1/2}(\Gamma_0)}^\vartheta.
\end{equation}
Apply the estimate to real and imaginary parts; no data on the rest of $\partial V$ are required. Using
\[
 \norm{\mathcal A_H^*r}_{\mathcal X_\Gamma}\leq C\norm r_{L^2(D_1)},
\]
equations \eqref{eq:dualglobal} and \eqref{eq:traceinterpolation} prove \eqref{eq:dualholder} with $\theta=\vartheta/2$.

The adjoint coefficients have the same support and satisfy
\[
 (c-\diver b)-\tfrac12\diver(-b)=c-\tfrac12\diver b\geq-\beta.
\]
Their bounds are controlled by the original class, so the argument also applies to $H^*$.
\end{proof}

\subsection{Runge approximation with the measured input norm}
The following argument implements quantitative Runge approximation directly for the slab. Its duality principle is related to the bounded-domain construction in \cite{RS}; here a quadratic minimization is used with the $H^{3/2}$ control norm.

\begin{lemma}[Approximation of extendible local solutions]\label{lem:runge}
There are $C>0$ and $\mu>0$ such that, for every solution $h\in H^1(D_2)$ of $Hh=0$ in $D_2$ and every $0<\epsilon<1$, one can find $f\in\mathcal X_\Gamma$ satisfying
\begin{align}
 \norm{P_H f-h}_{L^2(D_1)}&\leq\epsilon\norm h_{H^1(D_2)},\label{eq:rungeL2}\\
 \norm{P_H f-h}_{H^1(D_0)}&\leq C\epsilon\norm h_{H^1(D_2)},\label{eq:rungeH1}\\
 \norm f_{H^{3/2}(\Gamma_\ell)}&\leq C\epsilon^{-\mu}\norm h_{H^1(D_2)}.\label{eq:rungecost}
\end{align}
The same conclusions hold for $H^*$. The constants are uniform over the coefficient class.
\end{lemma}
\begin{proof}
Choose $\chi\in C_c^\infty(D_2)$ equal to one in a neighborhood of $\overline{D_1}$. For arbitrary $r\in L^2(D_1)$ and $w=(H^*)^{-1}(\mathbf1_{D_1}r)$, integration by parts gives
\[
 (h,r)_{L^2(D_1)}=(\chi h,H^*w)_{L^2(\Om)}
 =(H(\chi h),w)_{L^2(\Om)}.
\]
Local elliptic regularity makes $\chi h\in H^2(\Om)\cap H_0^1(\Om)$, so this calculation is valid for the stated weak solution. As $Hh=0$ and $b=0$ on $\supp\nabla\chi$,
\[
 H(\chi h)=-2\nabla\chi\cdot\nabla h-(\Delta\chi)h.
\]
Its support is a compact set $S\Subset D_2\setminus\overline{D_1}$. Lemma~\ref{lem:dualholder} consequently implies
\begin{equation}\label{eq:dualinequality}
 |(h,r)_{L^2(D_1)}|\leq C\norm h_{H^1(D_2)}
 \norm r_{L^2(D_1)}^{1-\theta}
 \norm{\mathcal A_H^*r}_{\mathcal X_\Gamma}^{\theta}.
\end{equation}

For $t>0$, minimize the strictly convex coercive functional
\[
 J_t(f)=\tfrac12\norm{\mathcal A_H f-h}_{L^2(D_1)}^2
       +\tfrac t2\norm f_{\mathcal X_\Gamma}^2.
\]
Its unique minimizer $f_t$ exists by the Hilbert space normal equation. Set $r_t=h-\mathcal A_H f_t$. That equation and comparison with $J_t(0)$ give
\[
 \mathcal A_H^*r_t=t f_t,\qquad
 \sqrt t\norm{f_t}_{\mathcal X_\Gamma}\leq\norm h_{L^2(D_1)},
\]
and
\[
 \norm{r_t}_{L^2(D_1)}^2+t\norm{f_t}_{\mathcal X_\Gamma}^2
 =\Rea(h,r_t)_{L^2(D_1)}.
\]
Applying \eqref{eq:dualinequality} to $r_t$ gives
\[
 \norm{r_t}_{L^2(D_1)}^{1+\theta}
 \leq C\norm h_{H^1(D_2)}^{1+\theta}t^{\theta/2},
 \qquad
 \norm{r_t}_{L^2(D_1)}\leq C\norm h_{H^1(D_2)}t^{\theta/(2+2\theta)}.
\]
The estimate also holds when $r_t=0$. Choose a sufficiently large fixed $C_0$ and set
\[
 t=(\epsilon/C_0)^{2(1+\theta)/\theta},\qquad
 \mu=(1+\theta)/\theta.
\]
This proves \eqref{eq:rungeL2} and \eqref{eq:rungecost}. Since $H(P_Hf_t-h)=0$ in $D_1$, the interior Caccioppoli estimate gives
\[
 \norm{P_H f_t-h}_{H^1(D_0)}
 \leq C\norm{P_H f_t-h}_{L^2(D_1)}.
\]
This proves \eqref{eq:rungeH1}. Replacing $H$ by $H^*$ proves the last assertion by Lemma~\ref{lem:dualholder}.
\end{proof}

\subsection{From the partial slab map to a full interior map}
Let $\Lambda^D_{b,c}:H^{1/2}(\partial D)\to H^{-1/2}(\partial D)$ denote the full ordinary DN map in $D$. Zero extension gives
\[
 H_0^1(D)\hookrightarrow H_0^1(\Om),\qquad
 \Rea\int_D\big(|\nabla v|^2+b\cdot\nabla v\,\overline v+c|v|^2\big)\,\dd x
 \geq(1-\beta/\lambda_1)\norm{\nabla v}_{L^2(D)}^2.
\]
A bounded trace lifting and the form argument imply
\[
 \norm u_{H^1(D)}\leq C\norm{u|_{\partial D}}_{H^{1/2}(\partial D)},\qquad
 \norm{\Lambda^D_{b,c}}_{H^{1/2}\to H^{-1/2}}\leq C.
\]
The weak normal trace is ordinary because the coefficients vanish near $\partial D$. We take $\langle\cdot,\cdot\rangle_{\partial D}$ linear in its first argument and conjugate-linear in its second.

\begin{proposition}[Quantitative interior-map transfer]\label{prop:interiormap}
Let the coefficients satisfy the uniform direct assumptions and have support in $K\Subset\Om$. For arbitrary nonempty open patches $G\subset\Gamma_1$ and $\gamma\subset\Gamma_\sigma$, set
\[
 d_\partial=\norm{\N_{b_1,c_1}-\N_{b_2,c_2}}_{\mathcal L(X_G,L^2(\gamma))}.
\]
There exist $C>0$ and $\rho\in(0,1)$ such that
\begin{equation}\label{eq:interiormap}
 \norm{\Lambda^D_{b_1,c_1}-\Lambda^D_{b_2,c_2}}
 _{\mathcal L(H^{1/2}(\partial D),H^{-1/2}(\partial D))}
 \leq C d_\partial^\rho.
\end{equation}
\end{proposition}
\begin{proof}
Choose disks $\Gamma^{\rm i}\Subset G$ and $\Gamma^{\rm o}\Subset\gamma$. For arbitrary $\varphi,\psi\in H^{1/2}(\partial D)$ of norm at most one, let
\[
 H_1h_1=0,\quad h_1|_{\partial D}=\varphi,
 \qquad H_2^*h_2=0,\quad h_2|_{\partial D}=\psi.
\]
The uniform solution bounds give $\norm{h_j}_{H^1(D)}\leq C$. Green's identity gives
\begin{equation}\label{eq:localidentity}
 \langle(\Lambda^D_1-\Lambda^D_2)\varphi,\psi\rangle_{\partial D}
 =\int_{D_0}((b_1-b_2)\cdot\nabla h_1+(c_1-c_2)h_1)\overline{h_2}\,\dd x.
\end{equation}
Subtract the $H_2$ solution with trace $\varphi$ from $h_1$ and test against $h_2$ to obtain the identity. Lemma~\ref{lem:runge} supplies
\[
 U_1=P_{H_1}f,\qquad U_2=P_{H_2^*}g,\qquad
 f\in\mathcal X_{\Gamma^{\rm i}},\quad g\in\mathcal X_{\Gamma^{\rm o}},
\]
with the estimates
\[
 \norm{U_j-h_j}_{H^1(D_0)}\leq C\epsilon,
 \qquad \norm{U_j}_{H^1(D_0)}\leq C,
 \qquad \norm f_{H^{3/2}}+\norm g_{H^{3/2}}\leq C\epsilon^{-\mu}.
\]
If the two Runge exponents differ, take their maximum. The corresponding identity on the slab is
\begin{equation}\label{eq:globalidentity}
 \int_{D_0}((b_1-b_2)\cdot\nabla U_1+(c_1-c_2)U_1)\overline{U_2}\,\dd x
 =\int_\gamma(\N_1-\N_2)f\,\overline g\,\dd S.
\end{equation}
The same subtraction proves this energy-space identity; $g$ is an adjoint test trace. The approximation bounds give
\begin{align*}
 &\left|\int_{D_0}\big((b_1-b_2)\cdot\nabla h_1+(c_1-c_2)h_1\big)\overline{h_2}\,\dd x\right.\\
 &\hspace{12mm}\left.-\int_{D_0}\big((b_1-b_2)\cdot\nabla U_1+(c_1-c_2)U_1\big)\overline{U_2}\,\dd x\right|
 \leq C\epsilon.
\end{align*}
The measured map controls the remaining boundary term:
\[
 \left|\int_\gamma(\N_1-\N_2)f\,\overline g\,\dd S\right|
 \leq d_\partial\norm f_{H^{3/2}}\norm g_{L^2(\gamma)}
 \leq C d_\partial\epsilon^{-2\mu}.
\]
Taking the supremum over $\varphi,\psi$ yields
\[
 \norm{\Lambda^D_1-\Lambda^D_2}\leq C(\epsilon+d_\partial\epsilon^{-2\mu}).
\]
For $0<d_\partial<1$, choose $\epsilon=d_\partial^{1/(2\mu+1)}$; for $d_\partial=0$, let $\epsilon\downarrow0$. Uniform boundedness of the local maps handles $d_\partial\geq1$. This proves \eqref{eq:interiormap} with $\rho=(2\mu+1)^{-1}$.
\end{proof}

\subsection{The full-data magnetic stability theorem and its application}
We use Caro and Pohjola \cite[Theorems~1.1 and~1.2]{CP}, allowing complex magnetic potentials. The required specialization assumes
\[
 \supp A_j\subset K\Subset D,\qquad
 \norm{A_j}_{L^\infty}+\norm{A_j}_{H^{1/2}(\R^n)}+\norm{q_j}_{L^\infty(D)}\leq C.
\]
For the full magnetic Cauchy data sets $\mathcal C_j$ on $\partial D$, those theorems give
\begin{align}
 \norm{\dd A_1-\dd A_2}_{H^{-1}(\R^n)}
 &\leq C[\log(2+\mathfrak d^{-1})]^{-c_*/(4n)},\label{eq:CPmag}\\
 \norm{q_1-q_2}_{H^{-1}(\R^n)}
 &\leq C[\log(2+\mathfrak d^{-1})]^{-c_*/(8n^3)},\label{eq:CPelectric}
\end{align}
for small $\mathfrak d=\operatorname{dist}(\mathcal C_1,\mathcal C_2)$ and a universal $c_*\in(0,1)$. The parameters in the cited theorems are
\[
 \text{Besov exponent}=\tfrac12,\quad \text{summability index}=2,\quad
 \lambda=1,\quad\text{auxiliary exponent}=\tfrac1n.
\]
Here $B^{1/2}_{2,2}=H^{1/2}$ with equivalent norms. Decreasing $c_*$ makes the same constant valid in both estimates.

We recall the distance to specify the interface. With $\mathcal T=H^1(D)/H_0^1(D)$ equipped with its quotient norm,
\[
 \mathfrak d=\max_{j\ne k}\ 
 \sup_{\substack{(f_j,g_j)\in\mathcal C_j\\\norm{f_j}_{\mathcal T}=1}}
 \inf_{(f_k,g_k)\in\mathcal C_k}
 \big(\norm{f_j-f_k}_{\mathcal T}+\norm{g_j-g_k}_{\mathcal T^*}\big).
\]
For the smooth domain $D$, these trace norms are equivalent to $H^{1/2}(\partial D)$ and $H^{-1/2}(\partial D)$. When both Dirichlet problems are uniquely solvable, choosing $f_k=f_j$ in the infimum gives
\begin{equation}\label{eq:cauchydistance}
 \mathfrak d\leq C\norm{\Lambda^D_{A_1,q_1}-\Lambda^D_{A_2,q_2}}
 _{H^{1/2}(\partial D)\to H^{-1/2}(\partial D)}.
\end{equation}

\begin{proof}[Proof of Theorem~\ref{thm:logcoeff}]
Set $A_j=ib_j/2$ and $q_j=Q_{b_j,c_j}$. Direct expansion gives
\[
 (D+A_j)^2+q_j=-\Delta+b_j\cdot\nabla+c_j.
\]
Interior support identifies the magnetic and ordinary DN maps. Coercivity gives unique solvability, and the coefficient bounds imply
\[
 \norm{A_j}_{L^\infty}+\norm{A_j}_{H^{1/2}(\R^n)}+\norm{q_j}_{L^\infty}\leq C.
\]
Thus \eqref{eq:CPmag}--\eqref{eq:CPelectric} apply. Equation~\eqref{eq:cauchydistance} and Proposition~\ref{prop:interiormap} give
\[
 \mathfrak d\leq C\norm{\Lambda^D_{b_1,c_1}-\Lambda^D_{b_2,c_2}}
 \leq Cd_\partial^\rho.
\]
For sufficiently small $d_\partial>0$, we obtain
\begin{equation}\label{eq:weakelliptic}
 \norm{\dd b_1-\dd b_2}_{H^{-1}(\R^n)}+
 \norm{Q_{b_1,c_1}-Q_{b_2,c_2}}_{H^{-1}(\R^n)}
 \leq C[\log(2+d_\partial^{-1})]^{-\varkappa_n},
 \quad\varkappa_n=\frac{c_*}{8n^3}.
\end{equation}
The power $d_\partial^\rho$ preserves the logarithmic rate, since
\[
 \log(2+(Cd_\partial^\rho)^{-1})\geq c\log(2+d_\partial^{-1})
\]
for small $d_\partial>0$. Equality of the local Cauchy data handles zero error; the uniform coefficient bounds handle large errors.

The additional Sobolev hypotheses imply
\begin{equation}\label{eq:strongprior}
 \norm{\dd b_j}_{H^s(\R^n)}+\norm{Q_{b_j,c_j}}_{H^s(\R^n)}\leq C.
\end{equation}
Indeed, fractional difference quotients and the derivative bound give
\[
 \norm{|b_j|^2}_{H^s}\leq C\norm{b_j}_{L^\infty}\norm{b_j}_{H^s},\qquad
 \norm{\diver b_j}_{H^s}\leq C\norm{b_j}_{H^{1+s}}.
\]
Interpolation on $\R^n$ gives, componentwise,
\[
 \norm v_{L^2}\leq\norm v_{H^{-1}}^{s/(1+s)}\norm v_{H^s}^{1/(1+s)}.
\]
Applying this to \eqref{eq:weakelliptic} and \eqref{eq:strongprior} yields
\begin{equation}\label{eq:strongelliptic}
 \norm{\dd b_1-\dd b_2}_{L^2}+
 \norm{Q_{b_1,c_1}-Q_{b_2,c_2}}_{L^2}
 \leq C[\log(2+d_\partial^{-1})]^{-\varkappa_n s/(1+s)}.
\end{equation}
Theorem~\ref{thm:quantmain} and Proposition~\ref{prop:transfer} now prove \eqref{eq:logcoeff} with \eqref{eq:logexponent}.

If the divergences agree, zero extension gives
\[
 \norm{\nabla(b_1-b_2)}_{L^2(\R^n)}^2=\norm{\dd b_1-\dd b_2}_{L^2(\R^n)}^2,
 \qquad \norm{b_1-b_2}_{H^1(\R^n)}\leq C\norm{\dd b_1-\dd b_2}_{L^2(\R^n)},
\]
by the Fourier identity and Poincar\'e's inequality on a fixed ball containing $K$. The equality
\[
 c_1-c_2=Q_{b_1,c_1}-Q_{b_2,c_2}-\tfrac14(b_1+b_2)\cdot(b_1-b_2)
\]
then gives \eqref{eq:logfixeddiv}.
\end{proof}

\begin{proof}[Proof of Corollary~\ref{cor:holderlog}]
The zero extensions retain their H\"older bounds because the supports stay away from $\partial\Om$. For compactly supported $f\in C^{0,\sigma_0}(\R^n)$,
\[
 \iint_{\R^n\times\R^n}\frac{|f(x)-f(y)|^2}{|x-y|^{n+2s}}\,\dd x\,\dd y
 \leq C\norm f_{C^{0,\sigma_0}}^2
 \left(\int_0^1 r^{2\sigma_0-2s-1}\,\dd r+\int_1^\infty r^{-2s-1}\,\dd r\right).
\]
Both integrals are finite for $0<s<\sigma_0$. Applying this to $c_j$, $b_j$ and the first derivatives of $b_j$ gives
\[
 \norm{b_j}_{H^{1+s}(\R^n)}+\norm{c_j}_{H^s(\R^n)}\leq C.
\]
Theorem~\ref{thm:logcoeff} applies; the scalar case follows from $Q_{0,c}=c$.
\end{proof}

\section{Uniqueness from an encoded boundary input}\label{sec:uniqueness}
\subsection{Construction and separation of the responses}\label{sec:encoding}
Fix $G$ as in \eqref{eq:G}. Recall that $X_G$ is the closure of $C_c^\infty(G)$ in $H^{3/2}(\Gamma_1)$. Choose a countable real sequence $\{g_k\}_{k\geq1}\subset C_c^\infty(G)$ whose complex linear span is dense in $X_G$.
To obtain the sequence, use separability of this closure and approximate a countable dense subset by smooth real functions. Including real and imaginary parts gives the asserted real sequence.

Fix $0<a<T_*<T$ and define
\[
 a_k=a+(T_*-a)(1-2^{-(2k-2)}),\qquad
 b_k=a+(T_*-a)(1-2^{-(2k-1)}),\qquad J_k=(b_k,a_{k+1}).
\]
Each $J_k$ is nonempty, and $a_k<b_k<a_{k+1}<T_*$. Let
$\theta\in C^\infty(\R)$ satisfy $\theta=0$ on $(-\infty,0]$ and $\theta=1$ on $[1,\infty)$. Set
\[
 \eta_k(t)=\theta\left(\frac{t-a_k}{b_k-a_k}\right).
\]
For $k\geq1$ choose the positive numbers
\begin{equation}\label{eq:dk}
 d_k=\frac{2^{-k}}{1+\displaystyle\sum_{m=0}^k
       \norm{\eta_k}_{C^m([0,T])}\norm{g_k}_{H^{m+2}(\Gamma_1)}}
\end{equation}
and prescribe
\begin{equation}\label{eq:F}
 F(t,x')=\sum_{k=1}^\infty d_k\eta_k(t)g_k(x').
\end{equation}

\begin{lemma}\label{lem:input}
The series \eqref{eq:F} defines a real smooth input with all the regularity and support properties in Theorem~\ref{thm:main}. For $t<T_*$ only finitely many summands are nonzero. For $t\geq T_*$ the input is constant in time. There is no loss of smoothness at $T_*$.
\end{lemma}
\begin{proof}
Fix nonnegative integers $r,s$ and take $m\geq\max\{r,s\}$. By \eqref{eq:dk},
\[
 d_k\norm{\eta_k g_k}_{C^r([0,T];H^s)}\leq2^{-k}\qquad(k\geq m).
\]
The finite initial sum and the summable tail therefore give
\[
 \sum_{k=1}^\infty d_k\eta_k g_k\quad\text{convergent in }C^r([0,T];H^s)
 \quad\text{for every }r,s.
\]
Termwise differentiation and Sobolev embedding prove smoothness. Since every ramp is constant near $T_*$,
\[
 \partial_t^rF(T_*)=0\quad(r\geq1),\qquad
 F(t)=\sum_{k=1}^\infty d_kg_k\quad(t\geq T_*).
\]
The ordering of the activation times and the support of $g_k$ give
\[
 F|_{[0,a]}=0,\qquad \supp F(t,\cdot)\subset\overline G\Subset\gamma_{\rm in}.
\]
Only finitely many terms are active at each $t<T_*$.
\end{proof}

\begin{proposition}[Recovery of a stationary partial map]\label{prop:decode}
Let $u_j$ be driven by \eqref{eq:F} and have arbitrary initial states $h_j\in X$. If \eqref{eq:dataeq} holds, then
\begin{equation}\label{eq:maps}
 \N_{b_1,c_1}^{\gamma_{\rm out}}g=
 \N_{b_2,c_2}^{\gamma_{\rm out}}g\qquad(g\in X_G).
\end{equation}
\end{proposition}
\begin{proof}
Before time $a$, the observed difference comes only from the initial states:
\[
 Z(t)=C_1S_1(t)h_1-C_2S_2(t)h_2=0\qquad(0<t<a).
\]
Graph-valued semigroup analyticity and the identity theorem imply
\[
 Z(t)=0\qquad(t>0).
\]
For each $k$, let $v_{j,k}$ be the zero-initial-state solution with input $\eta_kg_k$, extended by its plateau, and set
\[
 Y_k(t)=\partial_\nu v_{1,k}(t)|_{\gamma_{\rm out}}
        -\partial_\nu v_{2,k}(t)|_{\gamma_{\rm out}}.
\]
Only the first input is active on $J_1$. Thus
\[
 0=Z(t)+d_1Y_1(t)=d_1Y_1(t)\quad(t\in J_1),\qquad
 Y_1(t)=0\quad(t>b_1),
\]
where the last equality follows by analyticity. If the assertion holds for indices less than $k$, causality gives
\[
 0=Z(t)+\sum_{j=1}^kd_jY_j(t)=d_kY_k(t)\qquad(t\in J_k).
\]
Since $d_k>0$, analyticity again yields
\[
 Y_k(t)=0\qquad(t>b_k).
\]
Induction and the stationary limit \eqref{eq:limit} show that
\[
 (\N_{b_1,c_1}^{\gamma_{\rm out}}-\N_{b_2,c_2}^{\gamma_{\rm out}})g_k=0
 \qquad(k\geq1).
\]
Boundedness and density extend the equality to $X_G$, proving \eqref{eq:maps}.
Each interval $J_k$ uses a finite sum; the full solution across $T_*$ follows from Lemma~\ref{lem:input} and Proposition~\ref{prop:wellposed}.
\end{proof}

\subsection{The elliptic inverse step}\label{sec:elliptic}
We state exactly the external result used here. Write $D=-i\nabla$ and use the bilinear square $A\cdot A=\sum A_j^2$ for complex magnetic potentials.

\begin{theorem}[Krupchyk--Lassas--Uhlmann, specialized form]\label{thm:KLU}
Let $A_j\in W^{1,\infty}(\Om;\C^n)$ and $q_j\in L^\infty(\Om;\C)$ have compact support in $B\cap\overline\Om$. Suppose energy zero is admissible for $(D+A_j)^2+q_j$ and its real transpose $(D-A_j)^2+q_j$. If their partial magnetic DN maps at zero agree, for Dirichlet data on an open subset of $\Gamma_1$ containing $B\cap\Gamma_1$ and observations on an open subset of $\Gamma_\sigma$ containing $B\cap\Gamma_\sigma$, then
$\dd A_1=\dd A_2$ and $q_1=q_2$.
\end{theorem}
This is \cite[Theorem~1.1]{KLU} for $\sigma=2$ and \cite[Theorem~1.2]{KLU} for $\sigma=1$, at
\[
 k=0,\qquad k^2<\pi^2/L^2.
\]
Below this threshold, invertibility supplies admissibility and an $H^2(\Om)$ solution. The proof below checks invertibility for both the operator and its real transpose, as required by elliptic Runge approximation.

\begin{proof}[Proof of Theorem~\ref{thm:main}]
Choose
\begin{equation}\label{eq:conversion}
 A_j=\frac{i}{2}b_j,\qquad q_j=Q_{b_j,c_j}.
\end{equation}
Expanding the magnetic differential expression gives
\begin{align*}
 (D+A_j)^2+q_j
 &=-\Delta-2iA_j\cdot\nabla-i\diver A_j+A_j\cdot A_j+q_j\\
 &=-\Delta+b_j\cdot\nabla+\tfrac12\diver b_j
       -\tfrac14|b_j|^2+Q_{b_j,c_j}\\
 &=-\Delta+b_j\cdot\nabla+c_j=H_{b_j,c_j}.
\end{align*}
The regularity and compact-support hypotheses of Theorem~\ref{thm:KLU} follow from \eqref{eq:reg} and \eqref{eq:supp}. Since $b_j,c_j$ are real, the real transpose is also the Hilbert adjoint differential expression
\[
 (D-A_j)^2+q_j=-\Delta-b_j\cdot\nabla+c_j-\diver b_j.
\]
Lemma~\ref{lem:operator} shows that both Dirichlet realizations are invertible. Their free exterior threshold is positive and $k=0$ lies below it. Thus the admissibility requirements of Theorem~\ref{thm:KLU} hold.

The magnetic normal trace is
\[
 (\partial_\nu+iA_j\cdot\nu)u=\partial_\nu u-\tfrac12\rho u.
\]
The common Dirichlet data and normal trace cancel the correction. Proposition~\ref{prop:decode} therefore gives
\[
 (\partial_\nu+iA_1\cdot\nu)P_1g
 =(\partial_\nu+iA_2\cdot\nu)P_2g
 \quad\text{on }\gamma_{\rm out},\qquad g\in X_G.
\]
Choose a relatively open $G_0$ with
\[
 \overline B\cap\Gamma_1\subset G_0,\qquad\overline{G_0}\subset G.
\]
Mollification in the flat coordinates shows that any $g\in H^{3/2}(\Gamma_1)$ supported in $G_0$ admits
\[
 g_m\in C_c^\infty(G),\qquad g_m\longrightarrow g
 \quad\text{in }H^{3/2}(\Gamma_1).
\]
Hence $g\in X_G$, and all inputs required by Theorem~\ref{thm:KLU} are available. That theorem yields
\[
 \dd A_1=\dd A_2,\qquad q_1=q_2.
\]
Substitution of \eqref{eq:conversion} proves \eqref{eq:conclusion}.
\end{proof}

\begin{proof}[Proof of Corollary~\ref{cor:drift}]
Extend $b_1-b_2$ by zero to $\R^n$. Interior support introduces no boundary distributions, so
\[
 \dd(b_1-b_2)=0,\qquad\diver(b_1-b_2)=0\quad\text{in }\mathcal D'(\R^n).
\]
The Fourier identities give
\[
 |\xi|^2|\widehat{b_1-b_2}(\xi)|^2
 =\sum_{j<k}|\xi_j\widehat{(b_1-b_2)}_k-\xi_k\widehat{(b_1-b_2)}_j|^2
   +|\xi\cdot\widehat{b_1-b_2}|^2=0.
\]
Plancherel implies $b_1=b_2$. Equality of $Q_{b_1,c_1}$ and $Q_{b_2,c_2}$ then gives $c_1=c_2$. The scalar corollary is the case $b_1=b_2=0$ of Theorem~\ref{thm:main}.
\end{proof}

\subsection{A compact-class modulus}
\begin{lemma}[Continuity of the stationary map]\label{lem:mapcontinuous}
On the uniformly bounded coercive class of Theorem~\ref{thm:quantmain},
\[
 \norm{\N_{b_1,c_1}-\N_{b_2,c_2}}_{\mathcal L(X_G,Y)}
 \leq C\big(\norm{b_1-b_2}_{L^\infty}+\norm{c_1-c_2}_{L^\infty}\big).
\]
\end{lemma}
\begin{proof}
For $v=P_1g-P_2g$ the boundary values vanish and
\[
 H_1v=(H_2-H_1)P_2g=(b_2-b_1)\cdot\nabla P_2g+(c_2-c_1)P_2g.
\]
The uniform bound $H_1^{-1}:L^2\to H^2$, the uniform bound for $P_2$, and the normal trace theorem give the assertion.
\end{proof}

\begin{proof}[Proof of Theorem~\ref{thm:compactmain}]
Uniform support and H\"older bounds allow Arzel\`a--Ascoli on a fixed compact neighborhood of $K$. Every sequence in $\mathcal A$ has a subsequence with
\[
 b_k\longrightarrow b\quad\text{in }C^1(\overline\Om),\qquad
 c_k\longrightarrow c\quad\text{in }C^0(\overline\Om).
\]
Support, normal trace and coercivity pass to the limit. The difference-quotient definition also preserves the H\"older bound, so
\[
 (b,c)\in\mathcal A,\qquad \mathcal A\text{ is compact in }C^1\times C^0.
\]
Define the modulus directly, with the supremum of an empty set equal to zero:
\[
 \omega(s)=\sup\left\{
 \begin{aligned}
 &\norm{\dd b_1-\dd b_2}_{L^2}+\norm{Q_{b_1,c_1}-Q_{b_2,c_2}}_{L^2}:\\[-2pt]
 &(b_j,c_j)\in\mathcal A,\quad
 \norm{\N_{b_1,c_1}-\N_{b_2,c_2}}_{\mathcal L(X_G,Y)}\leq s
 \end{aligned}\right\}.
\]
This function is bounded and nondecreasing. Theorem~\ref{thm:KLU} and the common normal trace give
\[
 \N_{b_1,c_1}=\N_{b_2,c_2}
 \ \Longrightarrow\ \dd b_1=\dd b_2,\quad Q_{b_1,c_1}=Q_{b_2,c_2},
\]
as in the proof of Theorem~\ref{thm:main}; hence $\omega(0)=0$.

If $\omega(s)\not\to0$, compactness supplies convergent pairs $(b_{j,k},c_{j,k})$ and an $\epsilon_0>0$ such that
\[
 \norm{\N_{b_{1,k},c_{1,k}}-\N_{b_{2,k},c_{2,k}}}\longrightarrow0,\qquad
 \norm{\dd b_{1,k}-\dd b_{2,k}}_{L^2}
 +\norm{Q_{b_{1,k},c_{1,k}}-Q_{b_{2,k},c_{2,k}}}_{L^2}\geq\epsilon_0.
\]
Lemma~\ref{lem:mapcontinuous} makes the limiting maps equal. Continuity of the invariants gives the contradiction
\[
 0=\norm{\dd b_1-\dd b_2}_{L^2}+\norm{Q_{b_1,c_1}-Q_{b_2,c_2}}_{L^2}
 \geq\epsilon_0.
\]
Thus $\omega(s)\to0$, and Theorem~\ref{thm:quantmain} proves \eqref{eq:compactmain}.

For the last assertion, the zero extension of $b_1-b_2$ satisfies the Fourier identity
\begin{equation}\label{eq:hodge}
 \norm{\nabla(b_1-b_2)}_{L^2(\R^n)}^2
 =\norm{\dd b_1-\dd b_2}_{L^2(\R^n)}^2
  +\norm{\diver(b_1-b_2)}_{L^2(\R^n)}^2.
\end{equation}
The divergence term vanishes. Poincar\'e's inequality on a fixed ball containing $K$ gives
\[
 \norm{b_1-b_2}_{H^1}\leq C\norm{\dd b_1-\dd b_2}_{L^2}.
\]
Finally,
\[
 c_1-c_2=Q_{b_1,c_1}-Q_{b_2,c_2}
 -\tfrac14(b_1+b_2)\cdot(b_1-b_2),
\]
and the uniform drift bound proves \eqref{eq:fixeddivstable}.
\end{proof}

\subsection{A magnetic version}
Consider instead
\begin{equation}\label{eq:magheat}
 \partial_tu+\big((D+A)^2+q\big)u=0,
\end{equation}
where $A\in W^{1,\infty}(\Om;\R^n)$, $q\in L^\infty(\Om;\R)$ have the same support geometry and $q\geq-\beta$. The natural observation is $(\partial_\nu+iA\cdot\nu)u$. This equation is a complex-valued magnetic heat equation; it must not be identified with the ordinary real temperature equation when $A\ne0$.

\begin{corollary}\label{cor:magnetic}
For \eqref{eq:magheat}, there is a coefficient-independent input of the form \eqref{eq:F} such that equality of the magnetic normal traces on $(0,T)\times\gamma_{\rm out}$, with arbitrary $L^2$ initial states, implies $\dd A_1=\dd A_2$ and $q_1=q_2$. If $A_j\cdot\nu=0$ on the boundary, ordinary Neumann traces suffice.
\end{corollary}
\begin{proof}
The diamagnetic inequality, obtained by differentiating $|v|$ and regularizing at its zeros, gives
\[
 |\nabla|v||\leq|(\nabla+iA)v|\quad\text{a.e.},\qquad v\in H_0^1(\Om).
\]
By \eqref{eq:poincare} and $q\geq-\beta$,
\[
 \int_\Om\big(|(\nabla+iA)v|^2+q|v|^2\big)\,\dd x
 \geq\norm{\nabla|v|}_{L^2}^2-\beta\norm v_X^2
 \geq(\lambda_1-\beta)\norm v_X^2.
\]
The form is closed and equivalent to the $H^1$ norm, since
\[
 \norm{\nabla v}_{L^2}\leq\norm{(\nabla+iA)v}_{L^2}
 +\norm A_\infty\norm v_X.
\]
Its operator is positive self-adjoint. The regularity argument of Lemma~\ref{lem:operator} gives
\[
 D\big((D+A)^2+q\big)=H^2(\Om)\cap H_0^1(\Om).
\]
Replacing $A$ by $-A$ gives the same positivity and verifies admissibility of the real transpose at zero.

The magnetic observation is bounded on $H^2$ and, on zero-Dirichlet transients, satisfies
\[
 (\partial_\nu+iA\cdot\nu)v=\partial_\nu v.
\]
The lifting, stationary-limit and encoding arguments therefore apply. Theorem~\ref{thm:KLU} yields
\[
 \dd A_1=\dd A_2,\qquad q_1=q_2.
\]
The input may be real, since complex linearity extends equality from the real family to its complex span.
\end{proof}

\end{document}